\documentclass[12pt]{amsart}%
\usepackage{graphicx}
\usepackage{amscd}
\usepackage{amsmath}
\usepackage{amsfonts}
\usepackage{amssymb}%
\providecommand{\U}[1]{\protect\rule{.1in}{.1in}}
\providecommand{\U}[1]{\protect\rule{.1in}{.1in}}
\providecommand{\U}[1]{\protect\rule{.1in}{.1in}}
\newtheorem{theorem}{Theorem}
\theoremstyle{plain}

\newtheorem{corollary}{Corollary}

\newtheorem{definition}{Definition}

\newtheorem{proposition}{Proposition}
\newtheorem{remark}{Remark}

\numberwithin{equation}{section}

\begin{document}
\title[Normal numbers from Steinhaus' viewpoint ]{Normal numbers from Steinhaus' viewpoint }
\author{Daniel Pellegrino}
\address[Daniel Pellegrino]{ Depto de Matem\'{a}tica- Universidade Federal da
Para\'{\i}ba- Jo\~{a}o Pessoa, PB, Cep 58.051--900, Brazil }
\email{dmpellegrino@gmail.com}

\begin{abstract}
In this paper we recall a non-standard construction of the Borel sigma-algebra
$\mathcal{B}$ in $[0,1]$ and construct a family of measures (in particular,
Lebesgue measure) in $\mathcal{B}$ by a completely non-topological method.
This approach, that goes back to Steinhaus, in 1923, is now used to introduce
natural generalizations of the concept of normal numbers and explore their
intrinsic probabilistic properties. We show that, in virtually all the cases,
almost all real number in $[0,1]$ is normal (with respect to this generalized
concept). This procedure highlights some apparently hidden but interesting
features of the Borel sigma-algebra and Lebesgue measure in $[0,1].$

\end{abstract}
\maketitle

\section{Introduction}

The Borel sigma-algebra in $[0,1]$ is, in general, defined as the
sigma-algebra generated by the (open) intervals in $[0,1]$. So, we have a
natural \textquotedblleft topological\textquotedblright\ component in the
Borel sigma-algebra (and Lebesgue measure) in $[0,1].$ However, as it will be
shown, a completely different (non-topological) approach will be extremely
useful to deal with an extension of the concept of normal numbers and their
probabilistic aspects. The proposed characterization of the Borel
sigma-algebra in $[0,1],$ which goes back to 1923, with Steinhaus
\cite{Stein}, in the beginning of the conception of the modern probability
theory, shows a natural way to consider \textquotedblleft
weighted\textquotedblright\ measures in $[0,1],$ and to define and discuss
normality of numbers with respect to these measures.

The paper is organized as follows: In Section 2 we recall some background
results concerning product measures, in Section 3 we characterize the Borel
sigma-algebra in $[0,1]$ as a product sigma-algebra and introduce a family of
\textquotedblleft weighted\textquotedblright\ measures in $[0,1]$ and, in the
last section, we apply the previous results to generalize the concept of
normal numbers and explore their probabilistic behavior. The main goal of this
note is to translate, to the modern notation and terminology, Steinhaus'
striking ideas concerning normal numbers and also to call attention to his
contributions to the birth of modern probability theory.

Throughout, $\mathbb{R}$ denotes the set of all real numbers, $\mathbb{N}$
represents the set of all natural numbers $\{1,2,...\}$ and $\mathcal{B}$ is
the Borel sigma-algebra in the closed interval $[0,1]$. If $\Omega$ is any
set, $2^{\Omega}$ denotes the set of all subsets of $\Omega.$ If
$\mathcal{A}\subset2^{\Omega}$, then $\sigma(\mathcal{A})$ represents the
sigma-algebra generated by $\mathcal{A}$, $\#A$ denotes the cardinality of $A$
and $A^{\mathbb{N}}$ represents $A\times A\times\cdots.$

\section{Preliminary results}

Let $\Omega\neq\phi$ be a denumerable (or finite) set, $\mathcal{A}=2^{\Omega
}$ and let $\rho$ be a probability measure in $\mathcal{A}$.

Let
\[
(\Omega^{\mathbb{N}},\otimes\mathcal{A},\otimes\varrho)
\]
be the product space and denote%
\[
\mu=\mu_{\rho}=\otimes\varrho\text{ and }\mathcal{D}=\otimes\mathcal{A}%
\text{.}%
\]

\begin{definition}
$w=(a_{i})_{i=1}^{\infty}\in\Omega^{\mathbb{N}}$ is simply normal, with
respect to $\mu$, if
\[
\underset{n\rightarrow\infty}{\lim}\frac{S(w,r,n)}{n}=\rho(\{r\})\text{
}\forall r\in\Omega,
\]
where $S(w,r,n)$ denotes the total of indexes $i,1\leq i\leq n$, such that
$a_{i}=r.$
\end{definition}

\begin{definition}
$w=(a_{i})_{i=1}^{\infty}\in\Omega^{\mathbb{N}}$ is normal, with respect to
$\mu$, if
\[
\underset{n\rightarrow\infty}{\lim}\frac{S(w,B_{k},n)}{n}=%
{\displaystyle\prod\limits_{j=1}^{k}}
\rho(\{b_{j}\})\text{ }\forall k\in\mathbb{N},
\]

for each word $B_{k}=b_{1}...b_{k}$ of $k$ elements from $\Omega$, where
\[
S(w,B_{k},n)=\#\left\{  i\in\{1,...,n\};a_{i+j-1}=b_{j}\text{ for every
}j=1,...,k\right\}  .
\]

\end{definition}

The next two propositions are standard applications of the Strong Law of Large
Numbers. We sketch the proofs for the sake of completeness.

\begin{proposition}
The measure of the set
\[
\left\{  w=(a_{i})_{i=1}^{\infty}\in\Omega^{\mathbb{N}};w\text{ is simply
normal}\right\}
\]
is $1.$
\end{proposition}

Proof. Let $\mathcal{B}_{\mathbb{R}}$ denote the Borel sigma-algebra on
$\mathbb{R}$ and $r\in\Omega.$ Define, for every $n\in\mathbb{N}$,%
\begin{align}
X_{n}^{r}  &  :(\Omega^{\mathbb{N}},\mathcal{D},\mu)\rightarrow(\mathbb{R}%
,\mathcal{B}_{\mathbb{R}})\label{dem1}\\
w  &  =(a_{i})_{i=1}^{\infty}\mapsto\left\{
\begin{array}
[c]{c}%
0,\text{ if }a_{n}\neq r\\
1,\text{ if }a_{n}=r
\end{array}
\right. \nonumber
\end{align}

It is easy to see that each $X_{n}^{r}$ is measurable and hence $(X_{n}%
^{r})_{n=1}^{\infty}$ is a sequence of real random variables. \ Moreover,
$(X_{n}^{r})_{n=1}^{\infty}$ is an independent, integrable and identically
distributed sequence.

Note that%
\[%
{\displaystyle\int\nolimits_{\Omega^{\mathbb{N}}}}
X_{n}^{r}d\mu=\mu\left(  \left(  X_{n}^{r}\right)  ^{-1}(\{1\})\right)
=\rho(\{r\}).
\]
Hence, a well-known result due to Kolmogorov asserts that the sequence
$(X_{n}^{r})_{n=1}^{\infty}$ satisfies the Strong Law of Large Numbers, and
so
\[
\lim_{n\rightarrow\infty}\frac{1}{n}%
{\displaystyle\sum\limits_{i=1}^{n}}
X_{n}^{r}(w)=\rho(\{r\})\text{ (}\mu\text{-a.e)}.
\]
Therefore%
\[
\underset{n\rightarrow\infty}{\lim}\frac{S(w,r,n)}{n}=\rho(\{r\})\text{ (}%
\mu\text{-a.e)}.
\]

Denoting
\[
M^{r}=\left\{  w\in\Omega^{\mathbb{N}};\underset{n\rightarrow\infty}{\lim
}\frac{S(w,r,n)}{n}=\rho(\{r\})\right\}  ,
\]
we have that the set composed by the simply normal sequences on $\Omega
^{\mathbb{N}}$ is precisely%
\[%
{\displaystyle\bigcap\limits_{r\in\Omega}}
M^{r}.
\]
Since $\mu(\Omega^{\mathbb{N}}\smallsetminus M^{r})=0$ and $\Omega$ is, at
most, denumerable, we have%
\[
\mu\left(
{\displaystyle\bigcap\limits_{r\in\Omega}}
M^{r}\right)  =1.
\]
$\square$

\begin{proposition}
\label{prop2}The measure of the set $\left\{  w=(a_{i})_{i=1}^{\infty}%
\in\Omega^{\mathbb{N}};\text{ }w\text{ is normal}\right\}  $ is $1.$
\end{proposition}

Proof. Let $r_{1}....r_{k}$ be a word with $k$ elements from $\Omega$, with
$k\in\mathbb{N}$. Using the notation from (\ref{dem1}), define, for every
$n\in\mathbb{N}$,%
\[
\left\{
\begin{array}
[c]{c}%
Y_{(1),n}^{r_{1}...r_{k}}(w)=X_{kn-(k-1)}^{r_{1}}(w).X_{kn-(k-2)}^{r_{2}%
}(w)...X_{kn}^{r_{k}}(w)\\
Y_{(2),n}^{r_{1}...r_{k}}(w)=X_{kn-(k-2)}^{r_{1}}(w).X_{kn-(k-3)}^{r_{2}%
}(w)...X_{kn+1}^{r_{k}}(w)\\
\vdots\\
Y_{(k),n}^{r_{1}...r_{k}}(w)=X_{kn}^{r_{1}}(w).X_{kn+1}^{r_{2}}%
(w)...X_{kn+(k-1)}^{r_{k}}(w).
\end{array}
\right.
\]

It is plain that, for every $j=1,...,k,$ $(Y_{(j),n}^{r_{1}...r_{k}}%
)_{j=1}^{\infty}$ are integrable, independent and identically distributed
sequences. We also have, for every $j=1,...,k$,%
\[%
{\displaystyle\int\nolimits_{\Omega^{\mathbb{N}}}}
Y_{(j),n}^{r_{1}...r_{k}}d\mu=\mu\left(  \left(  Y_{(j),n}^{r_{1}...r_{k}%
}\right)  ^{-1}(\{1\})\right)  =%
{\displaystyle\prod\limits_{j=1}^{k}}
\rho(\{r_{j}\}).
\]
Hence, for every $j=1,...,k,$ $(Y_{(j),n}^{r_{1}...r_{k}})_{j=1}^{\infty}$
satisfy the Strong Law of Large Numbers and
\[
\lim_{n\rightarrow\infty}\frac{1}{n}%
{\displaystyle\sum\limits_{i=1}^{n}}
Y_{(j),n}^{r_{1}...r_{k}}(w)=%
{\displaystyle\prod\limits_{j=1}^{k}}
\rho(\{r_{j}\})\text{ (}\mu\text{-a.e).}%
\]
We thus have%
\[
\underset{n\rightarrow\infty}{\lim}\frac{S(w,r_{1}...r_{k},n)}{n}=%
{\displaystyle\prod\limits_{j=1}^{k}}
\rho(\{r_{j}\})\text{ (}\mu\text{-a.e).}%
\]
Denoting
\[
M^{r_{1}...r_{k}}=\left\{  w\in\Omega^{\mathbb{N}};\underset{n\rightarrow
\infty}{\lim}\frac{S(w,r_{1}...r_{k},n)}{n}=%
{\displaystyle\prod\limits_{j=1}^{k}}
\rho(\{r_{j}\})\right\}  ,
\]
the set composed by the normal sequences on $\Omega^{\mathbb{N}}$ is precisely%
\[%
{\displaystyle\bigcap\limits_{k\in\mathbb{N}}}
\left(
{\displaystyle\bigcap\limits_{r_{1}...r_{k}\in\Omega^{k}}}
M^{r_{1}...r_{k}}\right)  .
\]
Since $\mu(\Omega^{\mathbb{N}}\smallsetminus M^{r_{1}...r_{k}})=0$ and
$\Omega$ is, at most, denumerable, we have%
\[
\mu\left(
{\displaystyle\bigcap\limits_{k\in\mathbb{N}}}
\left(
{\displaystyle\bigcap\limits_{r_{1}...r_{k}\in\Omega^{k}}}
M^{r_{1}...r_{k}}\right)  \right)  =1.
\]
$\square$

\section{A family of measures in $[0,1]$}

The present section can be regarded, in some sense, as a translation of
Steinhaus' ideas from the roots of the modern probability theory (see
\cite{Stein}).

Henceforth, $\Omega=\{0,...,9\}$, $\mathcal{A}=2^{\Omega}$ and $\rho$ is a
probability measure on $\mathcal{A}$. If $(\Omega^{\mathbb{N}},\mathcal{D}%
,\mu_{\rho})$ is the product space, then $\mathcal{D}=\sigma(\mathcal{C}),$ with%

\[
\mathcal{C}=\left\{  \Omega\times\cdots\times\Omega\times\overset
{\text{position }n}{\{a\}}\times\Omega\times\cdots;a\in\Omega\text{ and }%
n\in\mathbb{N}\right\}  .
\]

Consider the mapping%
\begin{align}
\Psi_{\mu_{\rho}}  &  :(\Omega^{\mathbb{N}},\mathcal{D},\mu_{\rho}%
)\rightarrow([0,1],\sigma(\Psi_{\mu_{\rho}}(\mathcal{C})))\label{fim2}\\
w  &  =(a_{j})_{j=1}^{\infty}\mapsto\underset{j=1}{\overset{\infty}{\sum}%
}a_{j}10^{-j}.\nonumber
\end{align}
Note that $\Psi_{\mu_{\rho}}$ is \textquotedblleft almost
injective\textquotedblright, in the sense that, for every $(a_{j}%
)_{j=1}^{\infty}\in\Omega^{\mathbb{N}},$
\begin{align*}
&  \#\left(  \left\{  (b_{j})_{j=1}^{\infty}\in\Omega^{\mathbb{N}}%
;(b_{j})_{j=1}^{\infty}\neq(a_{j})_{j=1}^{\infty}\text{ and }\Psi_{\mu_{\rho}%
}((b_{j})_{j=1}^{\infty})=\Psi_{\mu_{\rho}}((a_{j})_{j=1}^{\infty})\right\}
\right) \\
&  =0\text{ or }1.
\end{align*}
Besides,%
\begin{align*}
&  \#\left(  \left\{  (a_{j})_{j=1}^{\infty}\in\Omega^{\mathbb{N}}%
;\exists(b_{j})_{j=1}^{\infty}\neq(a_{j})_{j=1}^{\infty}\text{ with }\Psi
_{\mu_{\rho}}((b_{j})_{j=1}^{\infty})=\Psi_{\mu_{\rho}}((a_{j})_{j=1}^{\infty
})\right\}  \right) \\
&  =\#(\mathbb{N}).
\end{align*}

From now on, sometimes we will write $\mu$ in the place of $\mu_{\rho}$ and
$\Psi_{\mu_{\rho}}$ will be simply denoted by $\Psi_{\mu}.$

Note that $\Psi_{\mu}$ is a random variable. In fact, if $A\in\Psi_{\mu
}(\mathcal{C}),$ then
\[
A=\Psi_{\mu}(\Omega\times\cdots\times\Omega\times\{b\}\times\Omega\times
\cdots)
\]
for some $b\in\Omega$ and thus%
\[
\left(  \Psi_{\mu}\right)  ^{-1}(A)=(\Omega\times\cdots\times\Omega
\times\{b\}\times\Omega\times\cdots)\cup D,
\]
where $D\in\mathcal{D}$ is a denumerable set. So,%
\[
\left(  \Psi_{\mu}\right)  ^{-1}(A)\in\mathcal{D}%
\]
and we can conclude that $\Psi_{\mu}$ is a random variable.

Denote the distribution of $\Psi_{\mu}$ by $\lambda_{\mu}.$ The measures
$\lambda_{\mu}$ in $\sigma(\Psi_{\mu}(\mathcal{C}))$ present an interesting
behavior, since, in general, these measures are \textquotedblleft
weighted\textquotedblright, i.e., they \textquotedblleft
protect\textquotedblright\ some digits. For example, it is not hard to see
that if%
\[
\rho(\{9\})=3/10,
\]
then%
\[
\lambda_{\mu}([\frac{9}{10},1])=3/10.
\]

Next, we will prove the following results, that are probably folkloric, but,
as far as we know, are (at least) very difficult to be found in the literature:

\begin{itemize}
\item $\Psi_{\mu}(\mathcal{D})=\sigma(\Psi_{\mu}(\mathcal{C}))=\mathcal{B}$.

\item If $\rho(\{r\})=1/10$ for every $r\in\Omega$, then $\lambda_{\mu}$ is
precisely the Lebesgue measure, i.e., when the measure $\rho$ is
\textquotedblleft non-weighted\textquotedblright, $\lambda_{\mu}$ coincides
with the Lebesgue measure.
\end{itemize}

\begin{theorem}
$\sigma(\Psi_{\mu}(\mathcal{C}))=\mathcal{B}$
\end{theorem}

Proof. Recall that%
\[
\mathcal{B}=\sigma\left(  \left\{  \left[  a,b\right]  ;a<b\text{ and }%
a,b\in\lbrack0,1[\right\}  \right)  .
\]

The following notation will be convenient:%
\begin{align*}
\left\{  >a\right\}   &  =\left\{  x\in\Omega;x>a\right\} \\
\left\{  <a\right\}   &  =\left\{  x\in\Omega;x<a\right\} \\
\left\{  \geq a\right\}   &  =\left\{  x\in\Omega;x\geq a\right\} \\
\left\{  \leq a\right\}   &  =\left\{  x\in\Omega;x\leq a\right\} \\
\left\{  >a,<b\right\}   &  =\left\{  x\in\Omega;x>a\text{ and }x<b\right\}  .
\end{align*}
We also define $\left\{  >a,\leq b\right\}  $, $\left\{  \geq a,<b\right\}  $
and $\left\{  \geq a,\leq b\right\}  $ in a similar way.

If $a,b\in\lbrack0,1[$, we can always write (uniquely)
\[
a=\underset{j=1}{\overset{\infty}{\sum}}a_{j}10^{-j}\text{ and }%
b=\underset{j=1}{\overset{\infty}{\sum}}b_{j}10^{-j}%
\]
with $\#\{j;a_{j}<9\}=\#\{j;b_{j}<9\}=\#(\mathbb{N}).$ Note that%
\[
\lbrack a,b]=%
{\displaystyle\bigcap\limits_{n=n_{0}+1}^{\infty}}
\left[  \underset{j=1}{\overset{n}{\sum}}a_{j}10^{-j},\underset{j=1}%
{\overset{n}{\sum}}b_{j}10^{-j}+\underset{j=n+1}{\overset{\infty}{\sum}%
}9.10^{-j}\right]  ,
\]
where $n_{0}$ is the smallest index for which $a_{n_{0}}\neq b_{n_{0}}.$

Since, for each $n\geq n_{0},$%
\begin{align*}
&  \left[  \underset{j=1}{\overset{n}{\sum}}a_{j}10^{-j},\underset
{j=1}{\overset{n}{\sum}}b_{j}10^{-j}+\underset{j=n+1}{\overset{\infty}{\sum}%
}9.10^{-j}\right] \\
&  =\Psi\left(  \{a_{1}\}\times\{a_{2}\}\times\cdots\times\{a_{n_{0}%
-1}\}\times\{>a_{n_{0}},<b_{n_{0}}\}\times\Omega\times\Omega\times
\cdots\right) \\
&
{\displaystyle\bigcup}
\Psi\left(
{\displaystyle\bigcup_{k=n_{0}+1}^{n}}
\{a_{1}\}\times\{a_{2}\}\times\cdots\times\{a_{k-1}\}\times\{>a_{k}%
\}\times\Omega\times\cdots\right) \\
&
{\displaystyle\bigcup}
\Psi\left(  \{a_{1}\}\times\{a_{2}\}\times\cdots\times\{a_{n}\}\times
\Omega\times\Omega\times\cdots\right) \\
&
{\displaystyle\bigcup}
\Psi\left(
{\displaystyle\bigcup_{k=n_{0}+1}^{n}}
\{b_{1}\}\times\{b_{2}\}\times\cdots\times\{b_{k-1}\}\times\{<b_{k}%
\}\times\Omega\times\cdots\right) \\
&
{\displaystyle\bigcup}
\Psi\left(  \{b_{1}\}\times\{b_{2}\}\times\cdots\times\{b_{n}\}\times
\Omega\times\cdots\right)  ,
\end{align*}
we can easily conclude that
\[
\left[  \underset{j=1}{\overset{n}{\sum}}a_{j}10^{-j},\underset{j=1}%
{\overset{n}{\sum}}b_{j}10^{-j}+\underset{j=n+1}{\overset{\infty}{\sum}%
}9.10^{-j}\right]  \in\sigma(\Psi_{\mu}(\mathcal{C}))
\]
and hence%
\[
\mathcal{B}\subset\sigma(\Psi_{\mu}(\mathcal{C})).
\]
Now, we must show that $\sigma(\Psi_{\mu}(\mathcal{C}))\subset\mathcal{B}$.

It suffices to show that $\Psi_{\mu}(\mathcal{C})\subset\mathcal{B}$. Let
$A\in\Psi_{\mu}(\mathcal{C})$. Hence%
\begin{align*}
A  &  =\Psi(\Omega\times\cdots\times\Omega\times\overset{\text{position }%
n}{\{b\}}\times\Omega\times\cdots)\\
&  =%
{\displaystyle\bigcup\limits_{\substack{a_{j}\in\{0,....,9\}\\j=1,...,n-1}}}
\left[  \underset{j=1}{\overset{n-1}{\sum}}a_{j}10^{-j}+\frac{b}{10^{n}%
},\underset{j=1}{\overset{n-1}{\sum}}a_{j}10^{-j}+\frac{b}{10^{n}}%
+\underset{j=n+1}{\overset{\infty}{\sum}}9.10^{-j}\right]  \in\mathcal{B}%
\text{.}%
\end{align*}
$\square$

\begin{theorem}
$\Psi_{\mu}(\mathcal{D})=\mathcal{B}$.
\end{theorem}

Proof. Using the previous result, all we need to show is that $\sigma
(\Psi_{\mu}(\mathcal{C}))=\Psi_{\mu}(\mathcal{D}).$

The proof that $\Psi_{\mu}(\mathcal{D})$ is a sigma-algebra needs a little bit
of hardwork, but is standard and we omit.

It is plain that%
\[
\sigma(\Psi_{\mu}(\mathcal{C}))\subset\Psi_{\mu}(\mathcal{D}).
\]
We will prove the converse inclusion.

It is not difficult to show that
\[
\mathcal{R}=\left\{  A\in\mathcal{D};\Psi_{\mu}(A)\in\sigma(\Psi_{\mu
}(\mathcal{C}))\right\}
\]
is a sigma-algebra and, since $\mathcal{C}\subset\mathcal{R}$, we have%
\[
\mathcal{D}=\sigma(\mathcal{C})\subset\mathcal{R}\text{.}%
\]
From the definition of $\mathcal{R}$ we conclude that%
\[
\Psi_{\mu}(\mathcal{D})\subset\sigma(\Psi_{\mu}(\mathcal{C}))
\]
and the proof is done. $\square$

Finally we have:

\begin{theorem}
\label{2}If $\rho(\{a\})=1/10$ for every $a\in\Omega$, the distribution of the
random variable $\Psi_{\mu_{\rho}}$ is the Lebesgue measure.
\end{theorem}

Proof. Let $\overline{\mu}$ be the distribution of $\Psi_{\mu_{\rho}}$ and
$\lambda$ be the Lebesgue measure on the Borel sigma-algebra on $[0,1].$ In
order to show that $\overline{\mu}$ and $\lambda$ coincide, it suffices to
show that they coincide over the intervals of $[0,1].$ In fact, in this case
they will coincide in the algebra $\mathcal{U}$ composed by the finite union
of disjoint intervals and, by invoking Carath\'{e}odory Extension Theorem, the
measures $\overline{\mu}$ and $\lambda$ will coincide in $\mathcal{B}%
=\sigma(\mathcal{U}).$

Note that the set%
\[
J=\left\{  x\in\lbrack0,1];\exists m\in\mathbb{N}\text{ so that }x=%
{\displaystyle\sum\limits_{j=1}^{m}}
x_{j}.10^{-j},\text{ }0\leq x_{j}\leq9\text{ (}\forall j=1,...,m)\right\}
\]

is dense in $[0,1].$ So, we just need to show that $\overline{\mu}$ and
$\lambda$ coincide over the intervals $[a,b]\subset\lbrack0,1],$ with $a,b\in
J.$

Moreover, there is no loss of generality in dealing with intervals of the form%
\[
I=\left[  \underset{j=1}{\overset{n}{\sum}}a_{j}10^{-j},\underset
{j=1}{\overset{n}{\sum}}b_{j}10^{-j}\right]  .
\]
Let $n_{0}$ be the smallest index such that $a_{n_{0}}\neq b_{n_{0}}.$ Hence%
\[
\lambda(I)=\underset{j=n_{0}}{\overset{n}{\sum}}(b_{j}-a_{j})10^{-j}.
\]
Consider $A\in\mathcal{D}$ given by%
\begin{align*}
A  &  =\left(  \{a_{1}\}\times\{a_{2}\}\times\cdots\times\{a_{n_{0}-1}%
\}\times\{>a_{n_{0}},<b_{n_{0}}\}\times\Omega\times\cdots\right) \\
&
{\displaystyle\bigcup}
\left(
{\displaystyle\bigcup_{k=n_{0}+1}^{n}}
\{a_{1}\}\times\{a_{2}\}\times\cdots\times\{a_{k-1}\}\times\{>a_{k}%
\}\times\Omega\times\cdots\right) \\
&
{\displaystyle\bigcup}
\left(  \{a_{1}\}\times\{a_{2}\}\times\cdots\times\{a_{n}\}\times\Omega
\times\Omega\times\cdots\right) \\
&
{\displaystyle\bigcup}
\left(
{\displaystyle\bigcup_{k=n_{0}+1}^{n}}
\{b_{1}\}\times\{b_{2}\}\times\cdots\times\{b_{k-1}\}\times\{<b_{k}%
\}\times\Omega\times\cdots\right) \\
&
{\displaystyle\bigcup}
\left(  \{b_{1}\}\times\{b_{2}\}\times\cdots\times\{b_{n}\}\times
\{0\}\times\{0\}\times\cdots\right)  .
\end{align*}
Hence $\left(  \Psi_{\mu_{\rho}}\right)  ^{-1}(I)=A\cup D$ with $\mu_{\rho
}(D)=0$ and
\[
\mu_{\rho}(A)=\frac{1}{10^{n_{0}}}(b_{n_{0}}-a_{n_{0}}-1)+\underset{k=n_{0}%
+1}{\overset{n}{\sum}}\frac{1}{10^{k}}(9-a_{k})+\frac{1}{10^{n}}%
+\underset{k=n_{0}+1}{\overset{n}{\sum}}\frac{b_{k}}{10^{k}}+0,
\]
and straightforward calculations show that%
\[
\mu_{\rho}(A)=\lambda(I).
\]
We thus have%
\[
\overline{\mu}(I)=\mu_{\rho}\left(  \left(  \Psi_{\mu_{\rho}}\right)
^{-1}(I)\right)  =\mu_{\rho}\left(  A\cup D\right)  =\mu_{\rho}(A)=\lambda
(I).
\]
$\square$

\section{A more general approach to normal numbers}

The notion of normal numbers (with respect to Lebesgue measure) was introduced
by E. Borel \cite{Bor}, in 1909, and, since then, several interesting
questions on normal numbers have been investigated and various intriguing
problems remain open (for example, the normality of $\sqrt{2}$). The results
of the previous sections turns natural to consider the concept of normal
numbers with respect to other measures than the Lebesgue measure on $[0,1]$.
In this section, as an application of the previous results, we generalize the
concept of normal numbers and obtain the measure of the sets of normal numbers
(with this generalized concept). In particular, we give an alternative simple
proof (essentially due to Steinhaus \cite{Stein}) for the fact that almost all
real numbers in $[0,1]$ are normal, with respect to Lebesgue measure
(different proofs of this result can be found, for example, in \cite{Bor},
\cite{Har}, \cite{Kui}, \cite{Ni}).

If $\Omega=\{0,...,9\}$, $\rho$ is a probability measure on $\mathcal{A}%
=2^{\Omega}$ and $\lambda_{\mu_{\rho}}$ is the distribution of the random
variable $\Psi_{\mu_{\rho}}$ defined in (\ref{fim2}), a\ number $\eta
\in\lbrack0,\infty\lbrack,$ represented in the decimal scale by
\begin{equation}
\eta=\overset{}{[\eta]+%
{\displaystyle\sum\limits_{j=1}^{\infty}}
a_{j}10^{-j},a_{j}\in\{0,...,9\},\forall j\in\mathbb{N}}\text{,} \label{mm}%
\end{equation}
with $[\eta]=\sup\{r\in\mathbb{N};r\leq\eta\}$ and $\#\{n;a_{n}<9\}=\#\left(
\mathbb{N}\right)  ,$ is said to be \textbf{simply normal} \textbf{(with
respect to }$\lambda_{\mu_{\rho}}$\textbf{) }when%

\begin{equation}
\underset{n\rightarrow\infty}{\lim}\frac{S(\eta,r,n)}{n}=\rho(\{r\})\text{
}\forall r\in\{0,...,9\}, \label{primm}%
\end{equation}
where $S(\eta,r,n)$ denotes the total of indexes $i,1\leq i\leq n$ such that
$a_{i}=r.$ A number $\eta,$ as in (\ref{mm}), is said to be \textbf{normal
(with respect to }$\lambda_{\mu_{\rho}}$\textbf{) }if\ %

\begin{equation}
\underset{n\rightarrow\infty}{\lim}\frac{S(\eta,B_{k},n)}{n}=%
{\displaystyle\prod\limits_{j=1}^{k}}
\rho(\{b_{k}\})\text{ }\forall k\in\mathbb{N}. \label{pri}%
\end{equation}
for each word $B_{k}=b_{1}...b_{k}$ of $k$ digits, where
\[
S(\eta,B_{k},n)=\#\left\{  i\in\{1,...,n\};a_{i+j-1}=b_{j}\text{ for every
}j=1,...,k\right\}  .
\]

In particular, if $\rho(\{a\})=\frac{1}{10}$ for every $a\in\Omega,$ this
concept is precisely Borel's original concept of normal numbers, with respect
to Lebesgue measure.

\begin{remark}
This generalized concept arises some interesting situations. For example, for
\textquotedblleft degenerate\textquotedblright\ cases, in which $\rho
(\{a\})=1$ for some $a\in\{0,1,...,9\},$ normal numbers are very special
numbers, with a strong preference to the digit $\{a\}.$ For example, if
$\rho(\{a\})=1$ for some $a\in\{0,1,...,9\},$ then
\[
0,a0aa0aaa0aaaa0aaaaa0...
\]
is normal with respect to $\lambda_{\mu_{\rho}}$.
\end{remark}

The next result shows that, in virtually all cases, the sets $N_{\mu_{\rho}}$,
of normal numbers in $[0,1]$, with respect to $\lambda_{\mu_{\rho}}$, are so
that $\lambda_{\mu_{\rho}}(N_{\mu_{\rho}})=1$, but there is one situation in
which $\lambda_{\mu_{\rho}}(N_{\mu_{\rho}})=0.$

\begin{theorem}
The measure of the set of normal numbers in $[0,1]$, with respect to
$\lambda_{\mu_{\rho}}$ is:

(a) $0$, if $\rho(\{9\})=1.$

(b) $1$, if $\rho(\{9\})<1.$
\end{theorem}

Proof. In the following, $N_{\mu_{\rho}}$ denotes the set of all normal
numbers in $[0,1]$, with respect to $\lambda_{\mu_{\rho}}$ and $M_{\mu_{\rho}%
}$ represents the set of all normal sequences in $\Omega^{\mathbb{N}},$ with
respect to $\mu_{\rho}.$

(a) If $\rho(\{9\})=1,$ then $\mu_{\rho}\left(  \{9\}^{\mathbb{N}}\right)
=1.$ We have%

\[
\left(  \Psi_{\mu_{\rho}}\right)  ^{-1}(N_{\mu_{\rho}})=D_{1},
\]
with
\[
D_{1}=M_{\mu_{\rho}}\smallsetminus\left\{  (a_{j})_{j=1}^{\infty};\text{
}\exists N\in\mathbb{N}\text{ such that }a_{n}=9\text{ for every }n\geq
N\right\}  .
\]
Hence $N_{\mu_{\rho}}\in\mathcal{B}$ and, since
\[
D_{1}\cap\{9\}^{\mathbb{N}}=\phi,
\]
we have%
\[
\lambda_{\mu_{\rho}}(N_{\mu_{\rho}})=\mu_{\rho}\left(  D_{1}\right)  =0.
\]

(b) If $\rho(\{9\})<1,$ note that%
\[
\left(  \Psi_{\mu_{\rho}}\right)  ^{-1}(N_{\mu_{\rho}})=M_{\mu_{\rho}}\cup
D_{2},
\]
with $\mu_{\rho}\left(  D_{2}\right)  =0.$ Hence%
\[
N_{\mu_{\rho}}\in\Psi_{\mu_{\rho}}(\mathcal{D})=\mathcal{B}%
\]
and by invoking Proposition \ref{prop2} we conclude that
\[
\lambda_{\mu_{\rho}}(N_{\mu_{\rho}})=\mu_{\rho}(M_{\mu_{\rho}})=1.
\]
$\square$

\begin{corollary}
(Borel \cite{Bor}) The set of normal numbers in $[0,1]$, with respect to
Lebesgue measure, has measure $1.$
\end{corollary}

\textbf{Acknowledgement.} This paper improves the author%
\'{}%
s dissertation under supervision of Professor M\'{a}rio C. Matos. The author
acknowledges Professor Matos for the important help and advice.

\end{document}